\documentclass[preprint]{imsart}

\usepackage[numbers]{natbib}
\usepackage[OT1]{fontenc}
\usepackage{amssymb, amsmath, amsfonts,amsthm}
\usepackage[colorlinks]{hyperref}
\usepackage{graphicx}

\usepackage{amsfonts}
\usepackage{epstopdf}
\usepackage{subfigure}
\usepackage{amssymb}
\usepackage{amsmath}
\usepackage{url}

\arxiv{math.PR/0000000}

\startlocaldefs
\numberwithin{equation}{section}
\theoremstyle{plain}

\newtheorem{proposition}{Proposition}[section]

\newtheorem{lemma}{Lemma}[section]
\newtheorem{lemmaa}{Lemma A.\hspace{-1ex}}
\newtheorem{remark}{Remark}[section]
\endlocaldefs

\def\num#1{\hbox{(\ref{#1})}}
\def\d{\displaystyle }

\def\Be{\begin{equation}}
\def\Eeq{\end{equation}}
\def\be{\begin{eqnarray}}
\def\eeq{\end{eqnarray}}
\def\ben{\begin{enumerate}}
\def\een{\end{enumerate}}
\def\ba{\begin{array}}
\def\ea{\end{array}}
\def\bt{\begin{tabular}}
\def\et{\end{tabular}}
\def\bc{\begin{center}}
\def\ec{\end{center}}
\def\bi{\begin{itemize}}
\def\ei{\end{itemize}}
\def\bd{\begin{document}}
\def\ed{\end{document}}
\def\nn{\nonumber\\ }
\def\nnn{\nonumber\eeq }
\def\lf{\lefteqn}
\DeclareMathOperator*{\argmax}{arg\,max}
\DeclareMathOperator*{\argmin}{arg\,min}
\def\choose#1#2{{\;\big( \ba{c} \vspace*{-1mm} #1 \\ \vspace*{1mm} #2 \ea\big)\;}}
\def\as{\,\mbox{ as }\,}
\def\If{\,\mbox{ if }\,}
\def\for{\,\mbox{ for }\,}
\def\and{\,\mbox{ and }\,}
\def\And{\,\mbox{ and }\,}
\def\textas{\,\mbox{ as }\,}
\def\textif{\,\mbox{ if }\,}
\def\textfor{\,\mbox{ for }\,}
\def\textand{\,\mbox{ and }\,}
\def\P{\mathbb P}
\def\zero{0\!\!\!\;\!\! 0}
\def\X{\mbox{\boldmath $X$}}
\def\Y{\mbox{\boldmath $Y$}}
\def\E{\mathbb E}
\def\R{\mathbb R}
\def\N{\mathbb N}
\def\Z{\mathbb Z}
\def\var{{\,{\mathrm{Var}}}}
\def\VAR{{\,\mathrm{VAR}}}
\def\cov{{\,{\mathrm{Cov}}}}
\def\std{{\,\mathrm{Std}}}
\def\corr{{\,\mathrm{Corr}}}
\def\x{\mbox{\boldmath $x$}}
\def\ep{\varepsilon}
\def\til{\!\!\!_{_\sim}}
\def\haha{$\stackrel{..}{\smile}$}
\def\l{\left\{ }
\def\r{\right\} }
\def\ll{\left.}
\def\rr{\right.}
\def\ra{\rightarrow}
\def\supp{\mbox{supp}}
\def\c{ \left|\ba{c} \vspace*{-5mm}\\ \ea\right.\!\!\!\!\!\! }
\def\cc{ \,\left|\ba{c} \vspace*{-5mm}\\ \ea\right.\!\!\!\!\! }
\def\*{$\cdot$}
\def\Xbar{\bar{X}}

\begin{document}
\def\N{\mathbb{N}}
\def\di{\displaystyle}
\def\indi{1\hspace{-1,1mm}{\rm I}}
\def\Fb{{\mathbb{F}}}
\def\Rb{{\mathbb{R}}}
\def\Tb{{\mathbb{T}}}
\def\Pr{{\mathbb{P}}}
\def\F{{\cal{F}}}
\def\zetab{\mbox{\boldmath$\zeta$}}
\def\Pc{{\cal{P}}}
\def\E{{\mathbb E}}
\def\Var{{\mathrm {Var}}}
\def\alph{a}
\def\bet{b}
\def\X{\mbox{\boldmath $X$}}
\def\Nb{{\mathbb N}}
\def\Zb{{\mathbb Z}}


\begin{frontmatter}
\title{{\large\bf
DETECTION AND ESTIMATION OF MULTIPLE TRANSIENT CHANGES
\thanksref{T1}}}
\runtitle{Detection and Estimation of Multiple Transient Changes}
\thankstext{T1}{Research of M. Baron is supported by the NSF grant 1737960 and
the DARPA grant HR0011-18-C-0051. Research of S. Malov is supported by
the RSF grant 20-14-00072}

\begin{aug}
\author{\fnms {Michael} \snm{Baron$^{1,}$}\thanksref{t1}\ead[label=e1]{baron@american.edu}}
\author{\fnms{Sergey V.} \snm{Malov$^{2,}$}\thanksref{t2}\ead[label=e2]{malovs@sm14820.spb.edu}}
\thankstext{t1}{E-mail: baron@american.edu }
\thankstext{t2}{E-mail: malovs@sm14820.spb.edu }

\runauthor{M. Baron $\&$ S.V. Malov }

\affiliation{American University, Washington D.C., USA}
\affiliation{St.-Petersburg State University, St.~Petersburg, Russia}

\address{
$^1$American University, Washington D.C., USA}
\address{
$^2$St.-Petersburg State University, St.~Petersburg, Russia}
\end{aug}

\begin{abstract}
 Change-point detection methods are proposed
 for the case of temporary failures, or transient changes, when an
 unexpected disorder is ultimately followed by a readjustment and
 return to the initial state. A base distribution of the
 ``in-control'' state changes to an ``out-of-control'' distribution for
 unknown periods of time. Likelihood based sequential and retrospective
 tools are proposed for the detection and estimation of each pair of
 change-points. The accuracy of the obtained change-point estimates is
 assessed. Proposed methods offer simultaneous control the familywise false
 alarm and false readjustment rates at the pre-chosen levels.
 \end{abstract}


\end{frontmatter}

\section{Introduction to transient changes}\label{sec:introduction}

Transient changes, or temporary disorders, refer to the situations
when an initial distribution of observed data changes to a different
one and eventually returns to the original state. The moments of change
are usually unexpected and a'priori unknown, the underlying distributions
may be known or unknown, but the ultimate return to the initial distribution
is assumed to be inevitable. In general, a data sequence may
experience one or more transient changes, which can be changes
in the mean value, variance, or other characteristics of the
observed process. This article focuses on the detection of
such changes and estimation of change-points.

There is a wide range of practical situations that are subject to
transient changes. Applications in signal and image processing
for the detection of finite signals
are mentioned in \cite{tartakovsky2021optimal}, with the detection
of space objects detailed in \cite{tartakovsky2021optimal}, Section 6.
Detection of transient changes appears useful in medical diagnostics
based on the heart rate variability \citep{bianchi1993time}.
Application to the monitoring of chemical concentrations in
drinking water is detailed in \cite{guepie2012sequential}, section 6.
Analysis of transient changes is important in industrial process
control and power systems, for the identification of in-control
and out-of-control periods; a specific application is described in
\cite{zhou2019improving}.

\begin{figure}\label{fig:spikes}
\begin{picture}(0,130)(0,0)
\put(-189,-70){\includegraphics{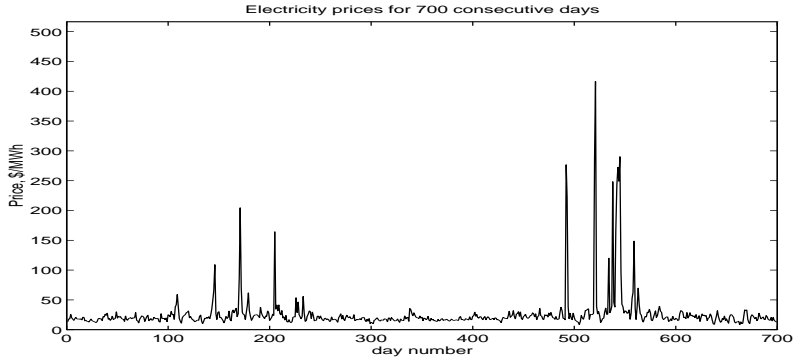}}
\end{picture}
\caption{Spikes in instantaneous electricity prices during two years
in the PJM (Pennsylvania--New Jersey--Maryland) energy market}
\end{figure}

Similar situations also occur in financial data from
deregulated energy markets. During the periods of high demand,
extreme weather conditions, maintenance or closure of a power plant,
the instantaneous price of electricity may experience a spike
lasting from several hours to several days, as on Figure~\ref{fig:spikes}. After each spike,
the distribution of prices returns to the initial state
\citep{BaronRosenbergSidorenko01,BaronRosenbergSidorenko02,tafakori2018forecasting,zhang2017regime}.
Accurate detection of spikes and estimation of their parameters
is needed for financial modeling and prediction that is critical
for proper valuation of energy options and contracts
\citep{RosenbergBryngelsonBaronPapalexopoulos10,RosenbergBryngelsonSidorenkoBaron02}.

Several statistical methods have been proposed for the detection
of transient changes. Under the assumption of a known
duration of the post-change period, the standard CUSUM algorithm for change-point
detection is modified and optimized in
\cite{guepie2017detecting,guepie2012sequential,tartakovsky2021optimal}.
The optimality is understood as the lowest probability of
missing a transient change \citep{guepie2017detecting,guepie2012sequential}
or the highest probability of detection \citep{tartakovsky2021optimal},
subject to the given probability
of a false alarm within the given time.
The optimized detection rule is the {\sl
window-limited CUSUM}, or WL-CUSUM. The special case of changes in the mean
is considered in \cite{noonan2020power},
where an approximate expression for the average run length to false
alarm is given for the {\sl moving-average sum} (MOSUM) algorithm.

When the assumption of a completely known duration of the period of change
is not realistic, one may consider it random, put a prior
distribution on each change-point, and consider the resulting
Bayesian problem, as in \cite{BaronRosenbergSidorenko02,
repin1991detection,tartakovskii1988detection}.

A more detailed overview of literature on
Bayesian and non-Bayesian transient change-point detection methods,
see \cite{egea2018performance,guepie2012sequential}.

In this work, we focus on the detection of transient changes and
estimation of change-points when the change intervals are completely
unknown. Algorithms are derived for the detection, estimation, and
testing of one transient change in Section 2, a known number of
transient changes in Section 3, and an unknown number of transient
changes in Section 4. The proposed detection method, a
self-correcting CUSUM procedure, is shown to detect transient
changes while controlling the familywise false alarm rate and the
familywise false readjustment rate simultaneously at the
pre-determined levels.

\section{Estimation and testing of one transient change interval}
\label{sec:single}

In this section, we assume at most one interval of change. Either
all the data follow the base distribution,
\[
H_0:\ X_1, \ldots,X_n\sim F,
\]
or there is one region of change $[\alph,\bet]$, so that
\[
H_1:\ \l\ba{lll}
 X_1,\ldots,X_\alph & \sim & F \\
 X_{\alph+1},\ldots,X_\bet & \sim & G \\
 X_{\bet+1},\ldots,X_n & \sim & F
 \ea\right.
\]
where $\alph$ and $\bet$ are unknown {\sl change-points}
while the distributions $F$ and $G$ are known.
The former case can be viewed
as the no-change null hypothesis $H_0$, and the latter as
the transient-change alternative $H_1$.

The goals are (1) to distinguish between $H_0$ and $H_1$
with a given level of significance, and
(2) to estimate change-point parameters $\alph$ and $\bet$
in the case of $H_1$.

\subsection{Maximum likelihood estimation}

For the parameter $(\alph,\bet)$, the log-likelihood
function is written as
 \be
  L(X;\alph,\bet)
&=& \sum_{i=1}^{\alph} \log f(X_i)+\sum_{i=\alph+1}^{\bet} \log
g(X_i) + \sum_{i=\bet+1}^{n} \log f(X_i) \nn
 &=&
\sum_{i=\alph+1}^{\bet} \log \frac{g(X_i)}{f(X_i)}
 + \sum_{i=1}^{n} \log f(X_i) \nn
 &=& S_b - S_a + \mbox{const},
 \label{LL-one-change}
 \eeq
 where $f$ and $g$ are probability densities of distributions
 $F$ and $G$ with respect to a reference measure $\mu$;
 \[
 S_t = \sum_{i=1}^{t} \log \frac{g(X_i)}{f(X_i)}
 \]
 is a random walk built on marginal log-likelihood ratios as
 its increments; and $\sum_{i=1}^{n} \log f(X_i)$ is a constant term
 as it does not depend on the unknown parameters $a$ and $b$.
 Measures $F$ and $G$ are not required to be mutually absolutely
 continuous, so that the log-likelihood ratio $\log(g/f)$
 assumes values in $\overline{\R}=[-\infty,\infty]$.

 Maximizing \num{LL-one-change}, we immediately obtain the
 maximum likelihood estimator (MLE)
 \Be\label{MLE-one-change}
 (\widehat\alph,\widehat\bet) = \argmax_{a\le b} (S_b-S_a).
 \Eeq

 According to \num{MLE-one-change}, the MLE returns the interval
 of the largest growth of random walk $S_t$. A direct method of calculating
 $\widehat\alph$ and $\widehat\bet$ can be proposed
 in terms of the associated cumulative-sum (CUSUM) process
 \Be\label{W}
 W_t=S_t-\min_{i\leq t} S_i,
 \Eeq
 which vanishes at every successive point of minimum of $S_t$.
 Given $\widehat\bet$, one finds $\widehat\alph$ by minimizing
 $S_t$ for $t\le \widehat\bet$, that is, finding the most recent
 zero of the CUSUM $W_t$. Then, the CUSUM does not return to zero
 between $\widehat\alph$ and $\widehat\bet$, and therefore,
 \Be\label{amplitude}
 S_{\widehat\bet}-S_{\widehat\alph}
  = W_{\widehat\bet}-W_{\widehat\alph}
  = W_{\widehat\bet}.
 \Eeq
 Maximizing \num{amplitude}, we obtain its
 computational formula for the MLE,
 \Be\label{MLE-W-one-change}
 \widehat\bet = \argmax W_t, \ \
 \widehat\alph = \max\l \mbox{Ker}(W) \cap [0,\widehat\bet) \r,
 \Eeq
 where
 Ker$(W)=\l t\,:\, W_t=0\r$ denotes the CUSUM's ``kernel'',
 or the set of its zeros.

\begin{figure}
 \scalebox{0.81}[0.77]{
 \begin{picture}(200,190)(60,-170)
 \put(-43,-200){\includegraphics{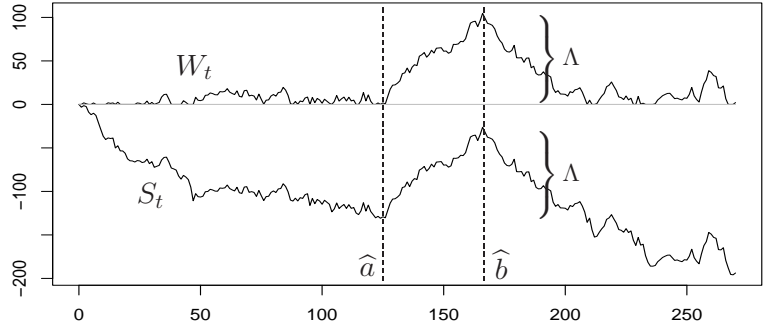}}
 \put(39,-111){\Large $S_t$}
 \put(57,-48){\Large $W_t$}
 \multiput(154,-149)(47,0){2}{\multiput(0,0)(0,4){34}{\line(0,1){2}}}
 \put(143,-147){\Large $\widehat\alph$}
 \put(206,-148){\Large $\widehat\bet$}
 \multiput(210,-43)(0,-57){2}{\large $\left.\begin{array}{c}
  \ \\[7mm] \
 \end{array}\right\}\Lambda$}
 \end{picture}}
 \caption{Maximum likelihood estimation of a single transient change
 interval. The likelihood-ratio test statistic $\Lambda$ is the largest increment
 of both processes $S_t$ and $W_t$.}
 \label{fig:S_and_W}
 \end{figure}

 As an illustration, an example of a log-likelihood ratio based
 random walk $S_t$, the associated CUSUM process $W_t$, and the
 resulting transient change-point estimator
 $(\widehat\alph,\widehat\bet)$ is shown in
 Figure~\ref{fig:S_and_W}.

 \subsection{Testing appearance of a transient change}

 The largest increment \num{amplitude} of the random walk $S_t$
 also serves as the log-likelihood ratio test statistic
 \[
 \Lambda = \log\frac{\max_{a<b}
 \l f(\X_{0:a})g(\X_{a,b})f(\X_{b:n})\r}{f(\X_{0:n})}
 = \max_{a<b} (S_b-S_a) = W_{\widehat\bet}
 \]
 for testing the no-change null  hypothesis against
 an alternative hypothesis that a transient change occurred
 in our data,
 \[
H_0:\ \X_{0:n}\sim F \ \mbox{ vs. } \
H_1:\ \l\ba{lll}
 \X_{0:\alph} & \sim & F \\
 \X_{\alph:\bet} & \sim & G \\
 \X_{\bet:n} & \sim & F
 \ea\right. \ \mbox{ for some } \alph < \bet
\]
 where $\X_{k:m} = (X_{k+1},\ldots,X_m)$ for any $k<m$.

 The likelihood-ratio test (LRT) rejects $H_0$ in favor of $H_0$
 if $\Lambda \ge h$ for some {\sl threshold
 $h$}. The choice of $h$ controls the balance between
 probabilities of Type I and Type II errors, or in other words,
 between the detection
 sensitivity and the rate of false alarms.

 In order to control the probability of a false alarm at the
 given level $\alpha$, we take advantage of the Doob's Maximal
 Inequality (for example, see \cite{Shiryaev95}, Section VII-3;
 \cite{stroock2013mathematics}, Section 7.1.1). It states that
 for a
submartingale $\{ Y_t\}$ and any constant $c\ge 0$,
\[
\Pr\left\{ \sup_{0\le t\le n} Y_t \ge c\right\} \le
\frac{\E(Y_n^+)}{c},
\]
where $x^+=\max\{ x,0 \}$.

The Doob's inequality can be applied directly to the LRT statistic
\[
\Lambda = W_{\widehat\bet} = \max_{0\le t\le n} W_t
\]
in the following way. The CUSUM process \num{W}
admits a recursive representation
\[
W_0 = 0, \ W_{t+1} = \max\l 0, \, W_t +
\log\frac{g(X_{t+1})}{f(X_{t+1})}\r
\]
(\cite{Page54}, Section 2.2).
Similarly, $U_t = \exp\l W_t\r$ can be expressed recursively as
$U_0 = 1, \ U_{t+1} = \max\{ 1, \, U_t g(X_{t+1}) / f(X_{t+1})\}$.
Therefore,
\[
\E_F\{ U_{t+1} \,|\, U_1,\ldots,U_t \} \ge U_t \E_F
\frac{g(X_{t+1})}{f(X_{t+1})} = U_t,
\]
showing that $U_t$ is a {\sl submartingale}.
Applying the Doob's maximal inequality to the process $\{
U_t\}$, we have
\be
 \Pr\l\mbox{ Type I error }\r
 &=& \Pr_F\l \Lambda \ge h\r
 = \Pr_F\{ \max_{0\le t\le n} W_t \ge h \}
 = \Pr_F\{ \max_{0\le t\le n} U_t \ge e^h \}
  \nn
 &\le& e^{-h}\E_F(U_n)
 = e^{-h}\E_F(e^{W_n}).
\label{Doob_bound_for_CUSUM}
 \eeq
 Thus, setting the threshold at
 \Be\label{h}
 h_\alpha = -\log\frac{\alpha}{E_F(e^{W_n})}
 \Eeq
 guarantees the probability of a false alarm no higher than $\alpha$.

 We conclude this section summarizing the obtained results.

 \begin{proposition}\label{prop:onechange}
 (Detection of a single transient segment)
 For the case of at most one transient change in the interval $[0,n]$,
 \ben\item
 The maximum likelihood estimator of the interval of change
 $[\alph, \bet]$ is given by \num{MLE-one-change} in terms of
 the random walk $S_t$ and by \num{MLE-W-one-change} in terms of the CUSUM
 process $W_t$.
 \item
 The likelihood-ratio test (LRT) rejects the no-change hypothesis if
 \[
 \Lambda = \max_{0\le t\le n} W_t \ge h
 \]
  for some threshold $h$.
 \item
 Threshold \num{h} yields a level $\alpha$ LRT.
 \item
 The change-point detection algorithm that reports a transient change
 at a stopping time $T_\alpha = \min\l t\,:\, W_t\ge h_\alpha\r$
 produces a false alarm with probability $\P\l\mbox{false alarm}\r\le\alpha$.
 \een
 \end{proposition}

 \subsection{Precision of the MLE}
\label{sec:MLEprop}

Introduce the direct and inverse shifted random walk processes $\{S_{\tau,i}^{+}\}_{i=1}^{n-\tau}$: $S_{\tau,i}^{+}=S_{\tau+i}-S_{\tau}$ and $\{S_{\tau,i}^{-}\}_{i=1}^{\tau-1}$: $S_{\tau,i}^{-}=S_{\tau-i}-S_{\tau}$ respectively; $S_0^+=S_0^-=0$, $S_i=S_i^+=S_{1,i}=$ and $S_i^-=S_{n,i}^-$, $i=1,\ldots,n$.
The MLE $(\hat\alph,\hat\bet)$ must satisfy inequalities
\begin{equation}
\label{equ:prelik}
\begin{array}{c}
\min\nolimits_{i}S_{\hat\alph,i}^{-}>0,
\qquad\min\nolimits_{i\leq \hat\bet-\hat\alph}S_{\hat\alph,i}^{+}> 0,\vspace{+.5em}\\
\max\nolimits_{i\leq \hat\bet-\hat\alph}S_{\hat\bet,i}^{-}< 0,\qquad
\max\nolimits_{i}S_{\hat\bet,i}^{+}<0.
\end{array}
\end{equation}
In the sequel, any pair $(\tilde\alph,\tilde\bet)$ that satisfies
(\ref{equ:prelik}) will be called a {\sl pre-likelihood estimator
(PLE)} in the sequel. In the next Section, we show that the PLE is
unique for sufficiently large $n$. If there exists a unique PLE, it
is equal to the MLE.

Let $\{W^-_i\}_{i=0}^{n}$: \ $W^-_0=0$, $W^-_{n-k}$
$=\max(0,S_1^-,\ldots,S_{k}^-)-S_{k}^-$, $k=1,\ldots,n$, be the
reverse CUSUM process; and $T=\{t: W_t=0\}$, $\widetilde
T=\{t:W^-_t=0\}$ and $\Tb=T\bigcup \widetilde T$ be the kernels of
the direct and reverse CUSUM $W_t$ and $W^-_t$. Note that
$W^-_i-W^-_j=W_j-W_i$ for any $i,j\in\{1,\ldots,n\}$. It is clear
that any PLE can be obtained from zeroes of the direct and inverse
CUSUM processes from the following characterization: a pair
$(\tilde\alph,\tilde\bet)$: $\tilde\alph<\tilde\bet$ is a PLE iff
$\tilde\alph\in T$, $\tilde\bet\in \widetilde T$ and $\{i\in
\Tb:\tilde\alph< i<\tilde\bet\}=\emptyset$.

In general, there may be multiple PLEs, and their number is random.
The exact distribution of the PLE is discussed in the Appendix.

We define a class of local estimators
$(\tilde\alph_{\gamma},\tilde\bet_{\gamma})$ that are constructed
around a given point $\gamma$, $$
 \tilde\alph_{\gamma}=\gamma-\arg\min\nolimits_{i\leq
\gamma} S_{\gamma,i}^-  \quad \mbox{and} \quad
\tilde\bet_{\gamma}=\gamma+\arg\max\nolimits_{i\leq n-\gamma}
S_{\gamma,i}^+.
$$
We call them local likelihood estimators (LLE). Every PLE coincides
with an LLE with respect to any point $\gamma$ inside the interval
defined by this PLE.

Below, we derive the distribution of an LLE with respect to a fixed
point $\gamma$.
%
%
Constructed this way, LLE $\tilde\alph_{\gamma}$ and
$\tilde\bet_{\gamma}$ are independent for any fixed $\gamma$. We
consider the interesting case of $\alph<\gamma\leq \bet$. The
distribution of $\tilde\bet_{\gamma}$ under $\bet-\gamma=k\geq 0$ is
the same as the distribution of MLE in the change point problem.
Then \arraycolsep=.4ex
\begin{eqnarray*}
\Pr(\tilde\bet_{\gamma}&=&\bet)=R_{G,k}^-(0)R_{F,n-\bet}^+(0)\\
&\geq & \exp\Bigl(-\sum\nolimits_{m=1}^{\infty} m^{-1}\Bigl( \Pr_F\bigl(\sum\nolimits_{j=1}^m Y_j \geq 0 \bigr)+\Pr_G\bigl(\sum\nolimits_{j=1}^m Y_j \leq 0 \bigr) \Bigr)\Bigr),
\end{eqnarray*}
where
 \Be\label{R_def}
 \ba{lll}
  R_{H,k}^+(x) &=& \Pr_H(\max(0,S_1^+,\ldots,S_k^+)<x),
 \\
  R_{H,k}^-(x) &=& \Pr_H(\max(0,S_1^-,\ldots,S_k^-)<x),
  \ea\Eeq
and the random walk process $S_k$ based on the sample from a distribution $H$.
The last inequality follows from the Spitzer's formula
(see \cite{Woodroofe82}) as $k,(n-b)\to\infty$ \cite{HuRukhin95}.
Moreover \cite{Hinkley70}, under $r>0$,
$$
\Pr(\tilde\bet_{\gamma}=\bet+r)=\int_0^{\infty}
R_{G,\bet-\gamma}^-(x)B^+_{F,r,n-b-r}(x)dx,
$$
and under $r<0$,
$$
\Pr(\tilde\bet_{\gamma}=\bet+r)= \int_0^{\infty}
R_{F,n-\bet}^+(x)B_{G,-r,\bet-\gamma+r}^-(x)dx,
$$
where
 \Be\label{B_def}\ba{lll}
B_{H,k,s}^+(y)dy &=& \Pr_H(\argmax\nolimits_{0\leq i\leq k+s}
S_i=k,S_k\in [y,y+dy))
 \\
B_{H,k,s}^-(y)dy &=& \Pr_H(\argmax\nolimits_{0\leq i\leq k+s}
(-S_i)=k,-S_k\in [y,y+dy)).
 \ea\Eeq
The distribution of $\tilde\alph_{\gamma}$ can be obtained in a
similar manner:
\begin{eqnarray*}
\Pr(\tilde\alph_{\gamma}&=&\alph)= R_{F,\alph}(0)R_{G,\gamma-\alph+1}^-(0)\\
&\geq & \exp\Bigl(-\sum\nolimits_{m=1}^{\infty} m^{-1}\Bigl( \Pr_G\bigl(\sum\nolimits_{j=1}^m Y_j \geq 0 \bigr)+\Pr_F\bigl(\sum\nolimits_{j=1}^m Y_j \leq 0 \bigr) \Bigr)\Bigr);
\end{eqnarray*}
$$
\Pr(\tilde\alph_{\gamma}=\alph+l)=\int_0^{\infty}
R_{G,\gamma-\alph+1}^-(x)B_{F,-l,\alph+l}^+(x)dx
$$
for $l<0$;  and
$$
\Pr(\tilde\alph_{\gamma}=\alph+l)=\int_0^{\infty}
R_{F,\alph}^+(x)B_{G,l,\gamma-\alph-l+1}^-(x)dx
$$
for $l>0$.

\begin{proposition}
\label{prp:1}
Let $\alph,\bet$ are fixed. Then
\arraycolsep=1pt
\begin{eqnarray*}
\sup_{\gamma\in (\alph,\bet]}\!\!\Pr(\tilde\bet_{\gamma}&=&\bet+r) \\
&\geq&
\begin{cases}
\exp\Bigl(\!-\!\!{\displaystyle\sum\limits_{m=1}^{\infty}} \frac{1}{m}\Bigl( \Pr_F\bigl({\displaystyle\sum\limits_{j=1}^m} Y_j \geq 0, \bigr)\!+\!\Pr_G\bigl({\displaystyle\sum\limits_{j=1}^m} Y_j \!\leq\! 0 \bigr) \Bigr)\Bigr)\;\;\mbox{for}\;\; r\!=\!0,\cr
{\displaystyle \int\nolimits_0^{\infty} R_{F,\infty}^+(x)B_{G,-r,\infty}^-(x)dx\quad\mbox{for}\;\; r<0,}\vspace{+.2em}\cr
{\displaystyle\int\nolimits_0^{\infty} R_{G,\infty}^-(x)B^+_{F,r,\infty}(x)dx\quad\mbox{for}\;\; r>0,}\cr
\end{cases}
\end{eqnarray*}
and
\begin{eqnarray*}
\sup_{\gamma\in (\alph,\bet]}\!\!\Pr(\tilde\alph_{\gamma}&=&\alph+l)\! \\
&\geq&
\begin{cases}
\exp\Bigl(\!-\!\!{\displaystyle\sum\limits_{m=1}^{\infty}} \frac{1}{m}\Bigl( \Pr_G\bigl({\displaystyle\sum\limits_{j=1}^m} Y_j \geq 0 \bigr)\!+\!\Pr_F\bigl({\displaystyle\sum\limits_{j=1}^m} Y_j \!\leq\! 0 \bigr) \Bigr)\Bigr)\;\;\mbox{for}\;\; l\!=\!0,\cr
{\displaystyle\int_0^{\infty} R_{G,\infty}^-(x)B_{F,-l,\infty}^+(x)dx\quad\mbox{for}\;\; l<0,}\vspace{+.2em}\cr
{\displaystyle\int_0^{\infty} R_{F,\infty}^+(x)B_{G,l,\infty}^-(x)dx\quad\mbox{for}\;\; l>0.}\cr
\end{cases}
\end{eqnarray*}

\end{proposition}

\begin{proof}
Continuing trajectories of the random walks  occurring at time $k$
in a neighborhood of $y$ we conclude that the probability of
reaching maximum at time $k$ is not decreased. Hence,
$B_{H,k,s}^{+}(x)$, $B_{H,k,s}^{-}(x)$ are non decreased in $s$ as
$s>0$. In a similar manner we obtain that $R_{H,s}^{+}(x)$ and
$R_{H,s}^{-}(x)$ are non decreased in $s$ as $s>0$.  The proposition
follows immediately.
\end{proof}

The right hand sides of inequalities in the last proposition for
$r\not=0$ or $s\not=0$ are quite hard for practical use. Its
approximations suitable for computation are obtained in
\cite{Hinkley70}.

The inequalities in Proposition \ref{prp:1} can be used immediately
to get the lower bound for the cumulative probabilities
$\Pr(s\leq\tilde\alph_{\gamma}-\alph\leq r)$ and
$\Pr(s\leq\tilde\bet_{\gamma}-\bet\leq r)$ for $s\leq r$ under
$\gamma\in (\alph,\bet]$. Bounds for the cumulative and tail
probabilities for LLE with respect to a detection point $\hat\gamma$
can be obtained from the following proposition.

\begin{proposition}
\label{prp:2} Let $\alpha\in (0,1)$; $\hat\gamma$ is the detection
point is such that $\Pr(\hat\gamma\in (\alph,\bet])\geq 1-\alpha$
for each $\alph<\bet$. Then
$$
\Pr(s\leq\tilde\alph_{\hat\gamma}-\alph\leq r)\geq
\sup\nolimits_{\gamma\in
(\alph,\bet]}\Pr(s\leq\tilde\alph_{\hat\gamma}-\alph\leq r)-\alpha
$$
$$
\Pr(s\leq\tilde\bet_{\hat\gamma}-\bet\leq r)\geq
\sup\nolimits_{\gamma\in
(\alph,\bet]}\Pr(s\leq\tilde\bet_{\hat\gamma}-\bet\leq r)-\alpha
$$
\end{proposition}

\begin{proof}
The proof is based on the Boole inequality
$$
1-\Pr(AB)=\Pr(\overline A \cup \overline B)\leq \Pr(\overline A)+\Pr(\overline B)
$$
that implies $\Pr(AB)\geq\Pr(A)-\Pr(\overline B)$. Then for any
fixed $a<b$,
\begin{eqnarray*}
\Pr(s\leq\tilde\alph_{\hat\gamma}-\alph\leq r)&\geq &\Pr(s\leq\tilde\alph_{\hat\gamma}-\alph\leq r, \hat\gamma\in (a,b]) \\
&\geq &\sup\nolimits_{\gamma\in (a,b]}\Pr(s\leq\tilde\alph_{\gamma}-\alph\leq r,\hat\gamma\in (a,b]) \\
&\geq &\sup\nolimits_{\gamma\in
(a,b]}\Pr(s\leq\tilde\alph_{\gamma}-\alph\leq r)-\alpha.
\end{eqnarray*}
The second inequality can be obtained analogously.
\end{proof}

\section{Asymptotic distribution of the MLE}


In this Section, we study the large-sample asymptotic behavior of
MLE $(\hat\alph,\hat\bet)$ as the sample size and all homogeneous
segments tend to infinity. We assume that the parameters
$\alph=\alph(n)$ and $b=\bet(n)$ and the interval of change $D =
[\alph(n),\bet(n)]$ are dependent on $n$, and
$\Delta=\min\l\alph(n),\bet(n)-\alph(n),n-\bet(n)\r\to\infty$ as
$n\to\infty$.  We use the notation $\Pr\equiv\Pr_{(D,n)}$ for the
distribution with a transient change and $\Pr_H$ for the case of
i.i.d. random variables $X_1,\ldots,X_n$ with the common
distribution function $H$. In particular,
$\Pr_F\equiv\Pr_{(\emptyset,n)}$ and $\Pr_G\equiv\Pr_{({\cal
N},n)}$, where ${\cal N}=\{1,\ldots,n\}$.

Next, we define random walks $S_k = \sum_{i=1}^k Y_i$ and
$\widetilde S_k=-\sum_{i=1}^k Y_i$. For example, for the transient
change-point detection problem, with log-likelihood ratios
$Y_i=\log\frac{g}{f}(X_i)$, the random walk $S_k$ is used to detect
a change from $F$ to $G$ whereas $\widetilde S_k$ is used to detect
a change from $G$ to $F$.

Let
$W_k$, $\widetilde W_k$ be the
corresponding CUSUM processes, where $W_0=\widetilde W_0=0$,
$W_k=S_k-\min_{i\leq k} S_k=(W_{k-1}+Y_k)\vee 0$, and $ \widetilde
W_k=\max_{i\leq k} S_k-S_k=(\widetilde W_{k-1}-Y_k)\vee 0 $.

We start with the following auxiliary results.

\begin{lemma}
\label{lem:1} Let $\E_F Y=c_1<0$.  Then for any $\epsilon >0$
$$
\Pr_F(\max\nolimits_{j\leq n} W_j/n>\epsilon)\to 0 \quad
\mbox{as}\quad n\to\infty.
$$
\end{lemma}

\noindent {\bf Proof.} Let ${\mathbb F}=\{{\cal F}_k\}_{k\in\Nb}$ be
the natural filtration associated with the process $Y_1,Y_2,\ldots$;
and $\tau_1,\tau_2,\ldots$ be the successive zeroes of the CUSUM
process $\{W_k\}_{k\in\Nb}$.

Introduce $Y^*_k=Y_k\indi_{\{Y_k\geq 0\}}$ and $S_k^*=\sum_{j=1}^k
Y^*_j$, $k\in\Nb$. Note that $\E Y^*_k=c_2<\infty$. By the Markov
property of the random process $\{S_k^*\}_{k\in \Nb}$:
$S_k^*=\sum_{j\leq k} Y_j^*$ with respect to the filtration
${\mathbb F}$ using Wald's identity, we obtain that $ \E_{F}
S^*_{\tau_k}=\E_{F} S^*_{\tau_1}=c_2 \E_{F} \tau_1 $ for all $k>1$.

Let $Z_1,Z_2,\ldots$ be independent copies of the random variable
$S^*_{\tau_1}$. Then
\begin{equation}
\label{equ:maxZ} \Pr_F(\max\nolimits_{j\leq n} W_j/n>\epsilon)\leq
\Pr_F(\max\nolimits_{j\leq n} Z_j/n>\epsilon).
\end{equation}
Note that $\E_{F} Z_1=c_3<\infty$, since $\E_{F} \tau_1<\infty$.
Denote, $H(u)=\Pr_F(Z_1\leq u)$ is the distribution function of
$Z$'s. It is well known that (see \cite{Leadbetter83}, sec. 1.5),
$$
n(1-H(u_n))\to\mu
$$
as $n\to\infty$, iff
$$
\Pr_F(\max\nolimits_{k\leq n}(Z_k)>u_n)\to 1- e^{-\mu}
\quad\mbox{as}\quad n\to\infty.
$$
Since $Z_1$ is non-negative, we obtain that $\E_F
Z_1=\int_0^{\infty} (1-H(x))dx<\infty$. Hence,
$$
\lim\sup_{t\to\infty} t(1-H(t))=\lim_{t\to\infty} t(1-H(t))=0,
$$
and, therefore, $u_n=o(n)$ as $n\to\infty$ for any $\mu>0$.
Therefore, the right hand side of the inequality (\ref{equ:maxZ})
is tended to $0$ as $n\to\infty$ for any $\epsilon>0$. The lemma is
proved.\quad\rule{2mm}{2mm}\medskip

\begin{lemma}
\label{lem:1as} Let $\E_F Y^2<\infty$. Then for any $\epsilon>0$,
$$
\Pr_F(\sup\nolimits_{m\geq n}\max\nolimits_{j\leq n}
W_j/n>\epsilon)\to 0 \quad \mbox{as}\quad n\to\infty.
$$
\end{lemma}

\noindent {\bf Proof.} Let $\epsilon>0$ be a fixed value. Then
$$
\Pr_F(\max\nolimits_{k\leq n} Z_k>n\epsilon)=1-H(n\epsilon)^n\sim
1-e^{-f(n)},
$$
and $f(n)=n(1-H(n\epsilon))=o(1)$, as $n\to\infty$. Hence,
$$
\Pr_F(\max\nolimits_{k\leq n} Z_k>n\epsilon)\sim f(n)
$$
as $n\to\infty$.  By the Borel--Cantelli lemma and Maclaurin--Cauchy
test, the convergence $\max\nolimits_{k\leq n} Z_k/n\to 0$ as
$n\to\infty$ holds  $\Pr_F$-almost sure if
$$
\int_{0}^{\infty} f(x)\,dx=\int_{0}^{\infty}
x(1-H(x\epsilon))\,dx<\infty.
$$
Note that  $E_F Y^2<\infty$ implies that $E_F Z_1^2<\infty$ in a
similar manner as in Lemma \ref{lem:1}. Then
$$
\E_{F} Z_1^2=2\int_{0}^{\infty} x(1-H(x))\,dx<\infty
$$
and, therefore,
$$
\int_{0}^{\infty}
x(1-H(x\epsilon))\,dx=\epsilon^{-2}\int_{0}^{\infty}
x(1-H(x))\,dx<\infty.
$$
Finally, we apply that
$$
\Pr_F(\sup\nolimits_{n\geq m}\max\nolimits_{j\leq n}
W_j/n>\epsilon)\leq \Pr_F(\sup\nolimits_{n\geq
m}\max\nolimits_{j\leq n} Z_j/n>\epsilon).
$$
The lemma is proved.\quad\rule{2mm}{2mm}\medskip

The next lemma follows immediately from the strong law of large
numbers (SLLN) and Lemma \ref{lem:1}.
\begin{lemma}
\label{lem:2} Let  $\E_F Y=c_1<0$, $\E_G Y=c_2>0$ and
$\lim\inf_{n\to\infty}\frac{\textstyle\Delta_n}{\textstyle n}\geq
\epsilon$ for some $\epsilon>0$.  Then
$$
\lim_{r\to\infty}\lim_{n\to\infty}\Pr(|\hat\alph_n-\alph|\geq r)=0;
\quad \lim_{r\to\infty}\lim_{n\to\infty}\Pr(|\hat \bet_n-\bet|\geq
r)=0.
$$
\end{lemma}

\noindent {\bf Proof.} For the most distant from $\alph$ version of
the PLE $\hat\alph_n$ we can write that
$$
\Pr(\alph-\hat\alph_n\geq r)=\Pr(\sup_{k\leq a-r}W_{k}=\sup_{k\leq
n}W_{k})\leq \Pr_F(\inf_{r\leq k\leq a}\widetilde S_k\leq
0)+\Pr(\sup_{k\leq a} W_k\geq S_{a,b-a})
$$
where
\begin{equation}
\label{equ:slln}
\begin{array}{c}
\di \Pr_F(\inf_{r\leq k\leq a}\widetilde S_k\leq
0)\leq\Pr_F(\inf_{k\geq r}\widetilde S_k\leq 0) \vspace{.5em}
\\
\di =\Pr_F(\sup_{k\geq r} S_k\geq 0)=\Pr_F(\sup_{k\geq r}
(S_k/k-c_1) \geq -c_1)\to 0
\end{array}
\end{equation}
as $r_0\to\infty$.  Moreover, for any $\delta>0$,
$$
\Pr(\sup\nolimits_{k\leq a} W_k\geq S_{a,b-a})\leq
\Pr_F(\sup\nolimits_{k\leq a} W_k>a\delta)+\Pr_G(S_{b-a}\leq
a\delta).
$$
The first term in the right-hand side of the last inequality tends
to $0$ as $n\to\infty$ by Lemma \ref{lem:1}, and the second term is
tended to $0$ as $n\to\infty$ by the law of large numbers since
$\lim\sup_{n\to\infty} a/(b-a)\leq\epsilon^{-1}$ and $c_2>0$.

Analogously, we obtain that $\Pr(\hat\bet_n-\bet\geq r) \to 0$ as
$r\to\infty$.

Let $\epsilon>0$; $r_{\epsilon}$ and $n_{\epsilon}$ are such that
$$
\Pr(\hat\bet_n-\bet> r_{\epsilon}) \leq\epsilon/2
$$
for all $n\geq n_{\epsilon}$. Since $\hat\alph_n\leq\hat\bet_n$, on
the event $A_{\epsilon}=\{\hat\bet_n-\bet\leq r_{\epsilon}\}$,
\begin{equation}
\label{equ:rsr} \Pr(\hat \alph_n-\alph\geq r)\leq\Pr_G(\sup_{r\leq
k\leq b-a} S_k\leq 0)+\Pr(S_{a,b-a}-\min_{k\leq r_{\epsilon}}
S_{b-a,k}\geq 0).
\end{equation}
The first term in the right hand side of the last inequality is
tended to $0$ as $r\to\infty$ uniformly in $n\geq 1$ as in
(\ref{equ:slln}). Then there exists an $r_{0\epsilon}$, such that
$\Pr_G(\sup_{r_{\epsilon}\leq k\leq b-a} S_k\leq 0)\leq \epsilon/4$.
Finally, $\min_{k\leq r_{\epsilon}} S_{b-a,k}=O_P(1)$, and,
therefore,
$$
\Pr(S_{a,b-a}-\min_{k\leq r_{\epsilon}} S_{b-a,k}\geq
0)=\Pr(S_{a,b-a}/(b-a)-c_2\geq -c_2+\min_{k\leq r_{\epsilon}}
S_{b-a,k}/(b-a))\to 0
$$
as $b-a\to\infty$. Hence, the second term in (\ref{equ:rsr}) is not
exceed $\epsilon/4$ under the sufficiently large $\Delta$. We
obtained that $ \Pr(\hat \alph_n-\alph\geq r_{\epsilon}\vee
r_{0\epsilon})\leq \epsilon, $ under the sufficiently large
$\Delta$, and, therefore,
$$
\lim_{r\to\infty}\lim_{\Delta\to\infty}\Pr(\hat \alph_n-\alph\geq
r)=0.
$$
Convergence $ \lim_{r\to\infty}\lim_{\Delta\to\infty}\Pr(\hat
\bet-\bet_n\geq r)=0 $ can be obtained in the similar manner.
Therefore, the lemma is proved.\quad\rule{2mm}{2mm}\medskip

Lemma~\ref{lem:2} yields the following proposition.

\begin{proposition}
\label{prp:lemp} Let $\gamma\!=\lambda \alph+(1-\lambda)\bet$,
$\hat\gamma\!=\lambda \hat \alph+(1-\lambda)\hat\bet$ for some
$\lambda\in (0,1)$; $\E_F Y<0$ and $\E_G Y>0$. Then
$$
\Pr_{\theta}(\hat \alph<\gamma<\hat\bet)\to 1 \quad \mbox{and}\quad
\Pr_{\theta}(\alph<\hat\gamma<\bet)\to 1,
$$
as $\lim\inf_{n\to\infty}\Delta_n/n\geq \epsilon$ for some
$\epsilon>0$. Moreover,
$$
\Pr_{\theta}(\max(\hat\bet-\gamma,\gamma-\hat\alph)>M)\to 1\quad
\mbox{and}\quad
\Pr_{\theta}(\max(\bet-\hat\gamma,\hat\gamma-\alph)>M)\to 1
$$
for any fixed $M>0$ for all $\theta=(\alph,\bet)$ as
$\lim\inf_{n\to\infty}\Delta_n/n\geq \epsilon$ for some
$\epsilon>0$.
\end{proposition}

\begin{remark}
(i). Lemma~\ref{lem:2} actually proves that the MLE
$(\hat\alph,\hat\bet)$ is unique for a sufficiently large $n$, and
if $\gamma\in (\alph,\bet)$ and  $(\tilde \alph_{\gamma},\tilde
\bet_{\gamma})$ is the LLE with respect to some point $\gamma$, then
$$
\Pr_{\theta}(\hat \alph=\tilde\alph_{\gamma},\hat\bet=\tilde
\bet_{\gamma}\;\mbox{under a sufficiently large}\;n)=1
$$
for all $\alph(n),\bet(n):\lim\inf_{n\to\infty}\Delta_n/n\geq \epsilon$ for some $\epsilon>0$.\\
(ii). Under some known point $\gamma$ between $\alph$ and $\bet$,
 the estimation problem reduces to two separate change-point estimation problems,
 on the direct $i=1,\ldots,\gamma$ and the inverse $i=n,n-1,\ldots,\gamma$ sets of indices. \\
(iii). The main results can be easily extended to the case of
multiple instability regions $D_n=\bigcup_{i=1}^{J_n} (a_j,b_j]$ as
$\Delta_n=\min_{j=0,\ldots,J_{n+1}}(b_j-a_j)\to\infty$ as
$n\to\infty$, where $a_0=0$ and $a_{J_n+1}=n$.
\medskip
\end{remark}


Remark 3.1(i) yields the following result, the asymptotic analogue
of Proposition 2.2, which establishes the asymptotic distribution of
the MLE $(\hat\alph,\hat\bet)$.
\begin{proposition}
\label{prp:1a}
Under the conditions of Lemma~\ref{lem:2}, for any fixed $r$ and
$l$, \arraycolsep=1pt
\begin{eqnarray*}
\lim_{n\to\infty}\Pr(\hat\bet&=&\bet+r)=p_r \\
&=&
\begin{cases}
\exp\Bigl(\!-\!\!{\displaystyle\sum\limits_{m=1}^{\infty}}
\frac{1}{m}\Bigl( \Pr_F\bigl({\displaystyle\sum\limits_{j=1}^m} Y_j
\geq 0, \bigr)\!+\!\Pr_G\bigl({\displaystyle\sum\limits_{j=1}^m} Y_j
\!\leq\! 0 \bigr) \Bigr)\Bigr)\;\;\mbox{for}\;\; r\!=\!0,\cr
{\displaystyle \int\nolimits_0^{\infty}
R_{F,\infty}^+(x)B_{G,-r,\infty}^-(x)dx\quad\mbox{for}\;\;
r<0,}\vspace{+.2em}\cr {\displaystyle\int\nolimits_0^{\infty}
R_{G,\infty}^-(x)B^+_{F,r,\infty}(x)dx\quad\mbox{for}\;\; r>0;}\cr
\end{cases}
\end{eqnarray*}
\begin{eqnarray*}
\lim_{n\to\infty}\Pr(\hat\alph &=&\alph+l)\!=\!q_l \\
&=&
\begin{cases}
\exp\Bigl(\!-\!\!{\displaystyle\sum\limits_{m=1}^{\infty}}
\frac{1}{m}\Bigl( \Pr_G\bigl({\displaystyle\sum\limits_{j=1}^m} Y_j
\geq 0 \bigr)\!+\!\Pr_F\bigl({\displaystyle\sum\limits_{j=1}^m} Y_j
\!\leq\! 0 \bigr) \Bigr)\Bigr)\;\;\mbox{for}\;\; l\!=\!0,\cr
{\displaystyle\int_0^{\infty}
R_{G,\infty}^-(x)B_{F,-l,\infty}^+(x)dx\quad\mbox{for}\;\;
l<0,}\vspace{+.2em}\cr {\displaystyle\int_0^{\infty}
R_{F,\infty}^+(x)B_{G,l,\infty}^-(x)dx\quad\mbox{for}\;\; l>0;}\cr
\end{cases}
\end{eqnarray*}
where $ \sum_{r\in\Zb} p_r=\sum_{l\in\Zb} q_l=1 $, and $R^+$, $R^-$,
$B^+$, and $B^-$ are defined by \num{R_def} and \num{B_def} in the
previous Section.
\end{proposition}


\section{Multiple transient changes and the familywise false alarm
rate}\label{sec:multiple}

Next, we consider a possibility of multiple transient changes
$[\alph_k,\bet_k]$, $k=1,\ldots,K$, where $K$ is
 the number of transient change intervals.
 The distribution of observed data oscillates between
 distributions $F$ and $G$, switching at unknown moments, so that
\[\begin{array}{lcccc}
 \X_{0:\alph_1} &\sim& X_1,\ldots,X_{\alph_1} &\sim& F \\
 \X_{\alph_1:\bet_1} &\sim& X_{\alph_1+1},\ldots,X_{\bet_1} &\sim& G \\
 \X_{\bet_1:\alph_2} &\sim& X_{\bet_1+1},\ldots,X_{\alph_2} &\sim& F \\
 \X_{\alph_2:\bet_2} &\sim& X_{\alph_2+1},\ldots,X_{\bet_2} &\sim& G \\
 \cdots &  & \cdots  &   & \cdots \\
 \X_{\bet_K:n} &\sim& X_{\bet_K+1},\ldots,X_{n} &\sim& F
 \end{array}
 \]
 One interpretation of this setting is a base distribution $F$,
 when the observed process is ``in control'', that is subject
 to sudden disorder times $\alph_k$, when it goes ``out of control'' to
 a disturbed distribution $G$. Each disorder will eventually be
 followed by an ``readjustment'' to the base distribution, which takes
 place at time $\bet_k$.

 The goal is to detect all the changes and estimate all
 $(2K)$ change-points $\alph_k$ and $\bet_k$. Facing a possibility
 of multiple changes, we aim to control a {\sl familywise false
 alarm rate} and a {\sl familywise false readjustment rate} that
 are understood as the probability of at least one erroneously
 detected change-point.

 That is, a $(2K)$-dimensional change-point parameter
 \[
  \l \alph_k, \bet_k \r_{k=1}^{k=K}
     = \l \alph_1, \bet_1; \ldots; \alph_K, \bet_K \r
 \]
 is estimated by a $2\widehat{K}$-dimensional estimator
 \[
  \l \widehat \alph_k, \widehat \bet_k \r_{k=1}^{k=\widehat K}
     = \l \widehat \alph_1, \widehat \bet_1; \ldots; \alph_{\widehat K}, \bet_{\widehat K}
     \r.
 \]
 A {\sl false alarm} is understood as an estimated segment
 $[\widehat\alph_k,\widehat\bet_k]$ that does not intersect
 with any disorder region $[\alph_m,\bet_m]$.
 The {\sl familywise false alarm rate} will be defined as the
 familywise error rate in the sense of
 \citep{HochbergTamhane87}, the probability of at least one
 false alarm,
 \Be\label{FAR}
 \mbox{FAR} = \Pr\l
   \cup_k \left( [\widehat\alph_k,\widehat\bet_k] \ \cap \
   ( \cup_j [\alph_j,\bet_j]) \ =\ \varnothing \right) \r.
 \Eeq
 Similarly, we call it a {\sl false readjustment}
   when the estimated ``in control'' interval
   $[\widehat\bet_k,\widehat\alph_{k+1}]$ does not contain
   any in-control observations, that is,
  \Be\label{FRR}
  \mbox{FRR} = \Pr\l
   \cup_k \left( [\widehat\bet_k,\widehat\alph_{k+1}] \ \cap\ \left(
   \cup_j [\bet_j,\alph_{j+1}]\right) \ =\ \varnothing \right)
   \r.
  \Eeq
  We aim at controlling the familywise rates of false alarms and false
  adjustments at pre-chosen levels $\alpha$ and $\beta$,
  respectively,
  \[
  \mbox{FAR} \le \alpha, \ \ \ \textand \ \ \ \mbox{FRR} \le \beta.
  \]
  We consider two situations, when the number of transient changes
  $K$ is known or unknown.

  \subsection{Known number of transient changes and MLE}\label{subsec:knownK}

The log-likelihood function of $\l
(\alph_k,\bet_k),\,k=1,\ldots,K\r$ is written as
 \[
 \ba{lll}\d L(X;\l (\alph_k,\bet_k)\r )
&=& \d\sum_{k=1}^{K} \sum_{i=\alph_k+1}^{\bet_k}
 \log\frac{g(X_i)}{f(X_i)}\ea
 \]

\begin{figure}
 \scalebox{0.88}[0.77]{
 \begin{picture}(200,160)(81,50)
 %
 \put(-32,16){\includegraphics{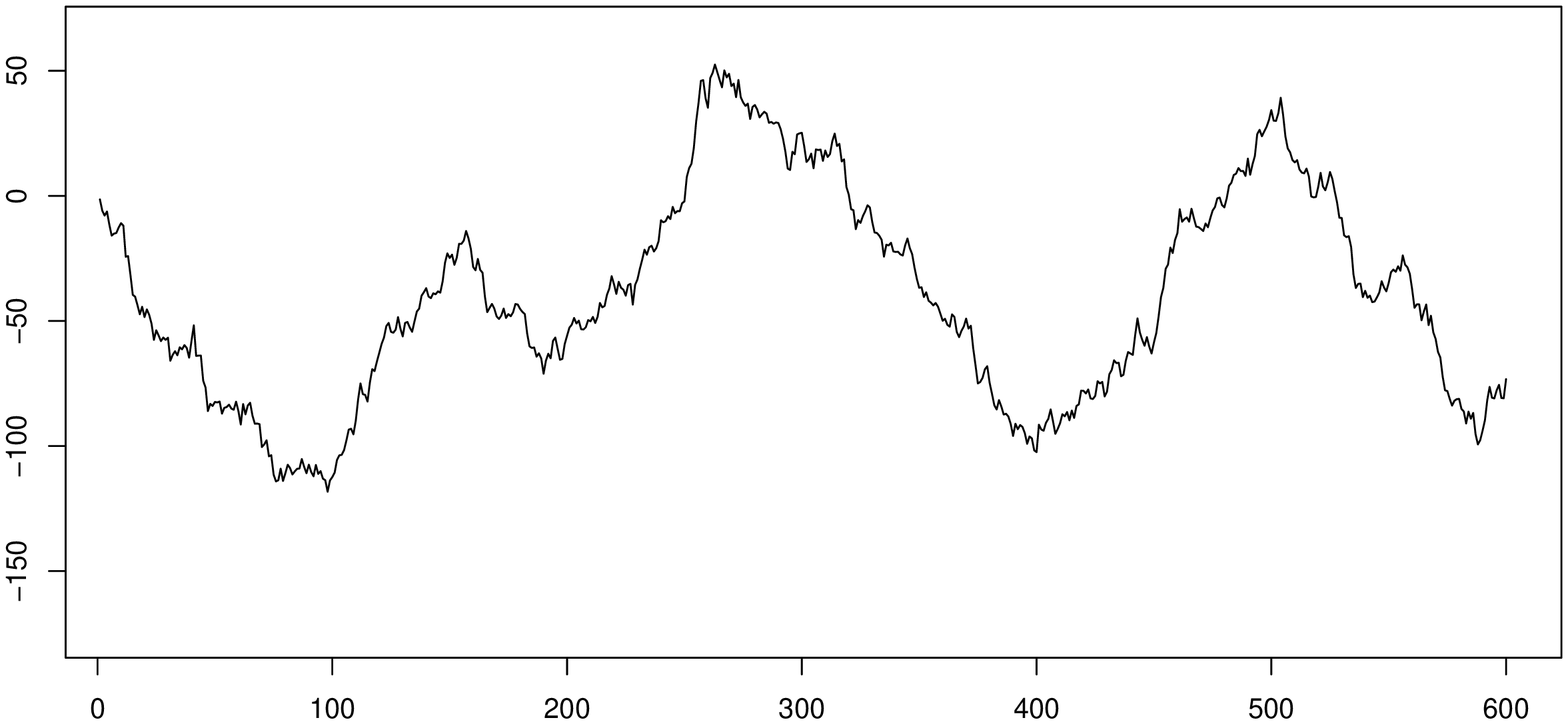}}
 \put(0,52.6){\vector(1,0){378}}
 \put(0,53){\vector(0,1){155}}
 \put(3,204){\Large $S_t$}
 \put(369,58){\Large $t$}
 \put(20,79){Step 1:}
 \put(20,68){Step 2:}
 \put(20,57){Step 3:}
 \put(60,79){$\widehat a_{1}$}
 \put(60,68){$\widehat a_{1}$}
 \put(60,57){$\widehat a_{1}$}
 \multiput(96,86)(0,4){15}{\line(0,1){2}}
 \put(92,79){---}
 \put(92,68){---}
 \put(93,57){$\widehat b_{1}$}
 \multiput(115,86)(0,4){9}{\line(0,1){2}}
 \put(111,79){---}
 \put(111,68){---}
 \put(113,57){$\widehat a_{3}$}
 \multiput(156,86)(0,4){25}{\line(0,1){2}}
 \put(150,79){$\widehat b_{1}$}
 \put(150,68){$\widehat b_{1}$}
 \put(150,57){$\widehat b_{3}$}
 \multiput(232,86)(0,4){5}{\line(0,1){2}}
 \put(228,79){---}
 \put(230,68){$\widehat a_{2}$}
 \put(230,57){$\widehat a_{2}$}
 \multiput(290,86)(0,4){22}{\line(0,1){2}}
 \put(285,79){---}
 \put(287,68){$\widehat b_{2}$}
 \put(287,57){$\widehat b_{2}$}
  \end{picture}}
 \caption{Estimation of multiple change-points}
 \label{fig:MLEmultiple}
 \end{figure}

 Maximizing it, we obtain the maximum likelihood estimator
 \[
 \l(\widehat\alph_k,\widehat\bet_k)\r_{k=1}^{k=K}
 = \argmax_{a_1 < b_1 < \ldots < a_k < b_K}
 \sum_{k=1}^K (S_{b_k}-S_{a_k}),
\]
 which are $K$ mutually disjoint
 intervals of the biggest growth of $S_t$
 (Figure~\ref{fig:MLEmultiple}).

A computational algorithm for $\l(\widehat\alph_k,\widehat\bet_k)\r$ can be obtained as an
iteration of steps outlined in Section~\ref{sec:single} for the
single-interval case, with a few modifications.

\emph{Step 1.} Apply \num{MLE-W-one-change} to obtain the first MLE
interval that corresponds to the interval of the biggest growth of
the random walk $S^{(1)}_t=S_t$ and the associated CUSUM process $W_t$,
\[
 \widehat\bet_1 = \argmax W_t, \ \
 \widehat\alph_1 = \max\l \mbox{Ker}(W) \cap [0,\widehat\bet) \r.
 \]

\emph{Step 2.} Apply Step 1 to the processes
\[
\ba{lllll}
S^{(2,1)}_t &=& S_t & \textfor & 1\le t\le \widehat\alph_1,
\\
S^{(2,2)}_t &=& -(S_t-S_{\widehat\alph_1}) & \textfor & \widehat\alph_1\le t\le \widehat\bet_1,
\\
S^{(2,3)}_t &=& S_t-S_{\widehat\bet_1} & \textfor & \widehat\bet_1 \le t\le n.
\ea\]
This results in three new intervals, $[c_1,d_1]$, $[c_2,d_2]$,
and $[c_3,d_3]$. Compare $\Delta_1 = S_{d_1}-S_{c_1}$,
$\Delta_2 = -(S_{d_2}-S_{c_2})$, and $\Delta_3 = S_{d_3}-S_{c_3}$,
and let $\Delta_j = \max\{\Delta_1,\Delta_2,\Delta_3\}$.

If $j=1$ or $j=3$,
add the corresponding interval to the MLE, i.e., let
\[
\widehat\alph_2 = c_j \textand \widehat\bet_2=d_j.
\]

If $j=2$, then let
\[
\widehat\bet_1 = c_2 \textand \widehat\alph_2=d_2,
\]
{\sl replacing} the previously found interval
$[\widehat\alph_1,\widehat\bet_1]$ with two intervals,
$[\widehat\alph_1,c_2]$ and $[d_2,\widehat\bet_1]$.

Based on the reversed log-likelihood ratios $\log(f/g)$,
the process $(-S_t)$ is actually the random walk
that can be used to detect a change from $G$ to $F$. Thus,
the found interval $[c_2,d_2]$ is a candidate for a
{\sl readjustment} period, a change back to the base distribution.
When the original walk $S_t$ drops more on $[c_2,d_2$ than
it grows on $[c_1,d_1]$ or $[c_3,d_3$, the sum of increments along
the obtained intervals
$[\widehat\alph_1,\widehat\bet_1]$ and $[\widehat\alph_2,\widehat\bet_2]$
is higher that on any other two intervals, and thus, they will
form the MLE for $K=2$.

\emph{Step $k$.} For $k=2,\ldots,K$, we repeat the same operations
as in Step 2. In every detected interval of change,
$[\widehat\alph_j,\widehat\bet_j]$, we find
an interval of the largest drop of $S_t$. In every interval
between them including $[1,\widehat\alph_1]$ and $[\widehat\bet_{k-1},n]$,
we find an interval of the largest growth of $S_t$.
Then we find the interval of the largest change among them.
If it is an interval of growth between
$\widehat\bet_j$ and $\widehat\alph_{j+1}$, we simply add
it to the list of intervals of change. If it is an interval $[c,d]$
of decrease between
$\widehat\alph_j$ and $\widehat\bet_j$, we replace
the previously found $[\widehat\alph_j,\widehat\bet_j]$
with two intervals, $[\widehat\alph_j,c]$ and $[d,\widehat\bet_j]$.

An example is shown in Figure~\ref{fig:MLEmultiple}.
At step 1, the interval of the largest growth is determined
as $[\widehat\alph_1=98,\ \widehat\bet_1=263]$.
At step 2, the second largest growth interval is determined as
$[\widehat\alph_2=400,\ \widehat\bet_1=504]$.
At step 3, the largest-growth interval with ends
at $c=157$ and $d=190$ is found inside
$[\widehat\alph_1, \widehat\bet_1]$. Therefore, we conclude
that a readjustment occurred between $c$ and $d$, and
$[\widehat\alph_1, \widehat\bet_1]=[98,263]$ is now replaced with
two intervals, $[98,157]$ and $[190,263]$.

\subsection{Unknown number of transient changes and
familywise error rates}

Since the number of changes $K$ is usually unknown, the algorithm in
Section~\ref{subsec:knownK} may either miss changes or produce
false alarms. As noted before, intervals of the biggest growth of
the random walk $S_t$ signal transient changes. Therefore, those
intervals where the increment in $S_t$ exceeds a certain threshold
will serve as estimated transient change intervals.

This threshold controls the rate of false alarms. As we show below,
no Bonferroni or Holm type correction is needed to control
the {\sl familywise} error rates. Instead, both the familywise
rate of false alarms \num{FAR} and the familywise rate of false
readjustments \num{FRR}
can be controlled by thresholds that are independent
of the true number of change-points, which can remain unknown.

The algorithm can be described as follows.

\bi\item[(i)] Introduce two CUSUM processes, renewed at a random
time $T\ge 0$,
\[
\ba{lllll} W_{T,t} &=& \d S_{T+t}-\min_{0\le i\le t}S_{T+i} &=&
\mbox{CUSUM based on $(S_{T+t}-S_T)$, renewed at $T$} \\
\widetilde W_{T,t} &=& \d \max_{0\le i\le t}S_{T+i}-S_{T+t} &=&
\mbox{CUSUM based on $-(S_{T+t}-S_T)$, renewed at $T$}
  \ea
 \]
 The CUSUM $W_{T,t}$ is set to detect the next disorder time,
 whereas $\widetilde W_{T,t}$ is tuned to determine the next
 readjustment time. A special case of $T=0$ results in
 the initial CUSUM processes $W_{t}$ and $\widetilde W_{t}$
 without any resetting.

\item[(ii)]
To control the familywise false alarm and false readjustment rates
at the desired levels $\alpha$ and $\beta$, respectively,
define thresholds as
 \Be\label{thresholds}
h_{\alpha} = -\log(\alpha \E_F^{-1}(e^{W_n}))
\ \textand \
\widetilde h_{\beta} = -\log(\beta \E_G^{-1}(e^{\widetilde W_n}))
 \Eeq

\item[(iii)]
The algorithm proceeds through the data series, detecting disorders
and readjustments at stopping times $\tau_k$ and post-estimating
change-points $\alph_k$ and $\bet_k$ sequentially for
$k=1,2,\ldots,K$ as follows,
 \[
 \ba{lll}
 \tau_1 &=& \inf\{t : 0<t\le n, W_t \ge h_{\alpha}\},
 \\
 \widehat\alph_1 &=& \max\l \mbox{Ker}\,W_t \cap [0,\tau_1)\r
 \\[3mm]
 \widetilde\tau_k &=& \tau_k + \inf\{t :
   0<t\le n-\tau_k,
   \widetilde W_{\tau_k,t} \ge \widetilde h_{\beta}\},
  \\
 \widehat\bet_k &=& \tau_k
   + \max\l\mbox{Ker}\,\widetilde W_{\tau_k,t} \cap [0,\widetilde\tau_k-\tau_k)\r
 ;
 \\[3mm]
 \tau_k &=&
   \widetilde\tau_{k-1} + \inf\{t :
   0<t\le n-\widetilde\tau_{k-1}, W_{\widetilde\tau_{k-1},t} \ge h_{\alpha}\},
  \\
 \widehat\alph_k &=& \widetilde\tau_{k-1}
   + \max\l\mbox{Ker}\, W_{\widetilde\tau_{k-1},t} \cap [0,\tau_k-\widetilde\tau_{k-1})\r
   ,
   \ea
  \]
  until $\tau_k=\infty$ or $\widetilde\tau_k=\infty$.
\ei
  By this definition of stopping times $\tau_k$, $\widetilde\tau_k$
  and change-point estimates $\widehat\alph_k$, $\widehat\bet_k$,
  each stopping time belongs to the corresponding interval of
  transient change that it is designed to detect,
  $\widehat\alph_k < \tau_k \le \widehat\bet_k$ and
  $\widehat\bet_{k-1} < \widetilde\tau_k \le \widehat\alph_k$.
  CUSUM processes $W_t$ and $\widetilde W_t$ are restarted and grounded
  at these times. As in the previous sections, change-points
  $\alph_k$ and $\bet_k$ are then estimated by
  the last zero points of
  restarted CUSUM processes
  $W_{\widetilde\tau_{k-1},t}$ and $\widetilde W_{\tau_k,t}$,
  respectively.

\begin{proposition}
\label{prop:FARcontrol} The transient change-point detection and
estimator scheme (i)-(iii) resulting in the estimator
$\{ \widehat\alph_k, \widehat\bet_k \}_{k=1}^{k=K}$ controls
familywise rates of false alarms and false readjustments at levels
\[
\mbox{FAR}\le\alpha \ \textand \ \mbox{FRR} \le \beta,
\]
for any unknown number of transient changes $K$.
\end{proposition}

\begin{proof}
According to the algorithm (i)-(iii), a false alarms occurs in the
interval $[\widehat\alph_k,\widehat\bet_k)$ if all the data in this
interval follow the distribution $F$, including the segment
$\X_{\widehat\alph_k:\tau_k}$ that triggered the false detection at
time $\tau_k$.

Also note that each renewed CUSUM process
 $W_{\tau_k,t} = S_{\tau_k+t}-\min_{0\le t\le \tau_k}S_{\tau_k+t}$ is dominated
 by the original CUSUM process $W_t$ on the corresponding segment,
 \[
 W_{\tau_k,t} \le W_{\tau_k+t}.
 \]
 This is because the subtracted term $\min_{t\ge 0}S_{t}$ in the
 original CUSUM process cannot exceed the corresponding minimum
 $\min_{0\le t\le \tau_k}S_{\tau_k+t}$ of the renewed CUSUM.

 Therefore, at least one false alarm can possibly occur only if the
 original CUSUM process $W_t$ exceeds the threshold $h_\alpha$ at
 least once in the interval $(0,n]$ under the distribution $F$.
 The probability of the latter event is bounded by the Doob's
 inequality. Similarly to \num{Doob_bound_for_CUSUM}, obtain
 \be
 \mbox{FAR} &\le& \P_F\l \bigcup_k \left( \max_{0<t\le
 (b_k-\widetilde\tau_{k-1})^+}W_{\tau_k,t} \ge h_\alpha\right)\r
 \le \P_F \l \max_{0<t\le n} W_t \ge h_\alpha\r
 \nn
 &\le& e^{-h}\E_F(e^{W_n})
 = \alpha,
 \nnn
 after substituting the first expression in \num{thresholds} for $h_\alpha$.

 The inequality $\mbox{FRR}\le \beta$ is proven along the same lines,
 replacing the CUSUM process $W_t$ with $\widetilde W_t$, and
 accordingly, the stopping times $\tau_k$ with $\widetilde\tau_k$
 and vice versa.

\end{proof}

\section{Nuisance parameters. Detection of changes to unknown distributions}

The case of known distributions $F$ and $G$ can sometimes apply to
real situations. For example, both distributions may be assumed to
belong to a parametric family $\l \F(\cdot | \theta),
\theta\in\Theta\r$ over a parameter set $\Theta$, where
$F(\cdot) = \F(\cdot | \theta_0)$ and
$G(\cdot) = \F(\cdot | \theta_1)$. The parameter $\theta_0$
of the base distribution $F$ may be known, representing the
``in-control'' state of the observed process. When the process
abruptly goes ``out of control'', the new distribution is
usually unknown. However, there is often the minimum magnitude of
a change $\Delta = |\theta_1-\theta_0|$ that is practically
reasonable to detect. The change-point detection algorithms
developed above can then be applied with known distributions
$\F(\cdot | \theta_0)$ and $\F(\cdot | \theta_1)$, where
$\theta_1 = \theta_0-\Delta$ or $\theta_1 = \theta_0+\Delta$ depending
on the direction of a change that is critical to be detected.
If both directions are important, two separate CUSUM procedures
or one two-sided CUSUM can be used \citep{crosier1986new}.

In this section, we elaborate the transient change-point detection
when changes of any magnitude are to be detected as long as the change
sustains for considerable time. At the same time, we achieve three
objectives:
\bi\item[--] The assumption of known distributions is lifted, they are allowed to be unknown;
\item[--] Distributions are no longer assumed the same during different transient change segments;
\item[--]Unknown nuisance parameters are estimated during the detection algorithm.
\ei

\subsection{The GLR process}

Suppose that the observed sequence follows a transient change-point model
with multiple segments from a family of distributions $\F(\cdot | \theta)$
with different (nuisance) parameters $\theta$. The baseline parameter $\theta_0$ may be
known or unknown; parameters of transient change segments $\theta_1, \theta_2,
\ldots$ are unknown, and they can be different for each transient change.

Following \cite{BaronRukhin97,reynolds2010evaluation,JJSiegmund92} and other work
on generalized likelihood ratio (GLR) processes,
the nuisance unknown parameters are to be estimated and replaced by the maximum
likelihood estimators. Due to unknown change-points, the segments of data
that should be used for estimating each nuisance parameter are also unknown.
Therefore, parameters are estimated for each {\sl potential change-point} $k$
and inserted into the equations for the random walk $S_t$, and subsequently,
the CUSUM $W_t$ for every $t$. The {\sl estimated CUSUM} process is then defined as
\Be\label{EstCUSUM}
 \widehat{W}_t = \max_{t^\omega \le k\le t-t^\omega}
   \sum_{i=k+1}^t \log\frac{f(\X_{k:t}|\widehat\theta_{k:t})}{f(\X_{0:k}|\widehat\theta_{0:k})}
   = \max_{t^\omega \le k\le t-t^\omega}
   \sum_{i=k+1}^t \log\frac{\d\max_{\theta\in\Theta} f(\X_{k:t}|\theta)}{\d\max_{\theta\in\Theta} f(\X_{0:k}|\theta)},
 \Eeq
 where $\widehat\theta_{0:k}$ and $\widehat\theta_{k:t}$ are maximum likelihood estimators
 of the pre- and post-change nuisance parameters, calculated from segments
 $\X_{0:k}$ and $\X_{k:t}$, as they would be estimated if a change occurred at time $k$.
 To avoid early false alarms caused by high variance of $\widehat\theta_{0:k}$ and $\widehat\theta_{k:t}$
 when they are estimated from short segments,
 the candidate change-point $k$ is separated from the ends of the interval $[0,t]$ by an amount
 $t^\omega$, where $\omega\in(0,1)$ will be determined later.

 A change-point is detected at the stopping time
 \[
 \widehat{T_h} = \inf\l t\le n | \widehat{W}_t \ge h\r.
 \]

 Of course, when the base distribution $\F(\cdot | \theta_0)$ is known, its known parameter $\theta_0$
 is used, and only $\theta_1$ is estimated.

 At the other end, when the family of distributions is unknown, and the problem is nonparametric,
 method similar to \num{EstCUSUM} can be used. As proposed in \cite{Baron00a}, the entire density
 is then treated as a nuisance parameter. Histogram density estimators can be used for the post-,
 and possibly, pre-change distributions, which actually converts the nonparametric problem into
 a parametric one, with multinomial distributions during each segment.

 \subsection{Probability of a false alarm}

 Can we still control the rate of false alarms when
 maximum likelihood estimators $\widehat\theta_0$ and
 $\widehat\theta_1$ are used in place of the exact nuisance
 parameters? Clearly, $\widehat\theta_j$ may differ substantially
 from $\theta_j$ for small $n$, resulting in substantial differences
 between random walks $S_t$ and $\widehat S_t$. However, one would
 expect the two random walks to be close for large $n$, when the
 nuisance parameters are estimated consistently.

 Suppose that the distributions $F$ and $G$ belong to a canonical
 exponential family with the common density
 \Be\label{expfamily}
 f(x|\theta) = e^{\theta x - \psi(\theta)} f(x|0).
 \Eeq
 Under this condition, Proposition 3.2 of \cite{BaronRukhin97}
 establishes the closeness of random walks $\hat S_t$ and $S_t$ that
 can be expressed as,
 \[
 \P\l \max_{k\le n^\gamma} |\hat S_k - S_k|\ge \epsilon\r = o\left(
 \exp\l -Cn^{(\omega-2\gamma)/3}\r\right),
 \]
 as $n\to\infty$, for some $C>0$ and any $\gamma \le \omega/2$.

 Applying this result to our change-point detection procedure with nuisance parameters, we
 deduce that for any $\epsilon>0$,
 \[
 \P\l \mbox{false alarm in } [0,n^\gamma]\r
 = \P_{\theta_0}\l \hat S_m - \hat S_k \ge h \mbox{ for some } k\le m\le
 n^\gamma\r
 \]\[
 \le
 \P_{\theta_0}\l S_m-S_k \ge h-2\epsilon \mbox{ or }
 |\hat S_m - S_m|\ge \epsilon \mbox{ or }
 |\hat S_k - S_k|\ge \epsilon  \mbox{ for some } k\le m\le
 n^\gamma\r
 \]\[
 \le o\left( 1 \right) + e^{-(h-2\epsilon)}\E_{\theta_0}{e^{W_n}}
 \le o\left( 1 \right) + \alpha e^{2\epsilon},
 \]
 as $n\to\infty$.

 The arbitrary choice of $\epsilon > 0$ yields that the probability
 of a false alarm is bounded by $\alpha + o(1)$. Further, for this
 probability bound, the power $\gamma$ can be taken to be
 $\omega/2$.

 \begin{proposition}
 When the underlying distributions $F$ and $G$ belong to an
 exponential family \num{expfamily} whose parameters are estimated
 from intervals of length at least $n^\omega$,
  \[
 \P\l \mbox{false alarm by GLR in } [0,\sqrt{n^\omega}]\r
 \le \alpha + o(1),
 \]
as $n\to\infty$.
 \end{proposition}


\appendix

\section{The distribution of PLE  single instability region}
\label{app}

Here we derive the exact distribution of PLE in the following
general setting.

Let $Y_1,Y_2,\ldots$ be a sequence of independent random variables
with the change points $0\leq a\leq \bet\leq\infty$; $Y_k$ has a
distribution function $G$ for $a<k\leq \bet$, and $F$ otherwise;
$S_k^+=\sum_{i=1}^k Y_i$ and $S_k^-=-\sum_{i=1}^k Y_i$ be the random
walks; $R_k^+=\max(0,S_1^+,\ldots,S_k^+)$ and
$R_k^-=\max(0,S_1^-,\ldots,S_k^-)$ be the sequential maxima;
$W_k^-=R_k^+-S_k^+$ and $W_k^+=R_k^--S_k^-=S_k-V_k$ be the CUSUM
processes; $V_k=\min(0,S_1^+,\ldots,S_k^+)$ be the sequential
minima, $k\in\N$. Introduce the following notations:
\begin{equation}
\label{equ:defR} R_{FG}(k,\bet)=\Pr(R_k^+=0) \quad{as}\quad  a=0
\end{equation}
and
\begin{equation}
\label{equ:defQ} Q_{FG}(k,\alph,\bet)=\Pr(V_k=0,W_k^-=0).
\end{equation}

The function $Q_{FG}(k,\alph,\bet)$ can be obtained recursively from
the following lemma.

\begin{lemmaa}
\label{lem:markov3}
Let $Y_1,\ldots,Y_n$ be independent random variables with the distribution functions $F_1,\ldots,F_n$ respectively;
$S_k=\sum_{j=1}^k Y_j$, \  \ $R_k=\max(0,S_1,\ldots,S_k)$, \  \ $V_k=\min(S_1,\ldots,S_k)$, \  \ $Q_k(v,w,s)=\Pr(V_k>v,S_k>s,R_k-S_k\leq w)$ \linebreak and $\di Q'_k(v,w,s)=\frac{\Pr(V_k>v,S_k\in [s,s+ds),R_k-S_k\leq w)}{ds}$, $k=1,\ldots,n$.
Then  $Q_1(v,w,s)=(F_1(w)-F_1(v\vee s))\indi_{\{w\geq v\vee s\vee 0\}},$
\begin{equation}
\label{equ:rec3}
Q_{k+1}(v,w,s)=\int Q_k(v,w+x,s\vee v-x)\,dF_{k+1}(x)\indi_{\{w\geq 0\}}
\end{equation}
and $Q'_1(v,w,s)=F([s,s+ds))\indi_{\{v\leq s, w\geq s\vee 0\}}$,
\begin{equation}
\label{equ:rec3d}
Q'_{k+1}(v,w,s)=\int Q'_k(v,w+x,s-x)\,dF_{k+1}(x)\indi_{\{ v\leq s,\;w\geq 0\}},
\end{equation}
$k=1,\ldots,n-1$.
\end{lemmaa}

\begin{proof}
The relation $Q_1(v,w,s)=(F_1(w)-F_1(v\vee s))\indi_{\{w\geq 0\}}$  is clear. Note that $V_{k+1}=\min(V_{k},S_k+Y_{k+1})$, $R_{k+1}=\max(R_k,S_k+Y_{k+1})$ and $S_{k+1}=S_k+Y_{k+1}$. Then $R_l-S_l-Y_{k+1}\leq t$ if $R_{k+1}=R_k$, and $R_k-S_k-Y_{k+1}\leq 0$ if $R_{k+1}=S_{k+1}$.
Therefore,
\arraycolsep=0.2ex
\begin{eqnarray*}
&&\Pr(V_{k+1}>v,S_{k+1}>s,R_{k+1}-S_{k+1}>w)
\\
&&=\Pr(V_k>v,S_k+Y_{k+1}>v,S_k+Y_{k+1}>s,R_k-S_k-Y_{k+1}\leq w)
\\
&&=\int \Pr(V_k>v,S_k>v\!\vee\! s-x, R_k-S_k\leq w-x)\indi_{\{w\geq 0\}} \,dF_{k+1}(x)
\end{eqnarray*}
for each $k\in 1,\ldots,n-1$. The recursive equations (\ref{equ:rec3d}) can be obtained in a similar manner.
The proof is completed.
\end{proof}

\begin{remark}
(i). Lemma A.\ref{lem:markov3} establishes the Markov property of the 3-dimensional process $(V_k,R_k,S_k)$, $k=1,\ldots,n$. \\
(ii). The recurrent formula (\ref{equ:rec3}) is applicable for fixed $v$, calculation of $Q_k(v,w,s)$ for all $v$ and each $k$ is not required in this case.
\end{remark}

Let $Y_1,\ldots,Y_n$ be a sequence with fixed change points $0\leq
a\leq b\leq\infty$; $(\tilde\alph_m,\tilde\bet_m)$, $m=1,\ldots,M$
be the PLEs. The total number $M$ of PLEs is random and the common
distribution of PLEs and $L$ is quite hard to be obtained. The
common distribution of a version $(\tilde\alph,\tilde\bet)$ of PLE
if $\Pr(M>1)$ is small can be evaluated by the following
probabilities obtained immediately from (\ref{equ:prelik}):
\arraycolsep=1pt
\begin{eqnarray*}
p_{lr}&=&\Pr(\cup_{m=1}^M\{\tilde\alph_m=\alph+l,\tilde\bet_m=\bet+r\})\\
&=& R_{FG}(\alph\!+\!l,l_*)Q_{FG}(d,l^*,d-r_*)R_{FG}(n-\bet-r,r^*)
\end{eqnarray*}
where $l_*=l\vee 0$, $l^*=-l\vee 0$, $r_*=r \vee 0$, $r^*=-r \vee
0$, $d=\bet-\alph+r-l$; $R_{FG}$ and $Q_{FG}$ are defined in
(\ref{equ:defR}) and (\ref{equ:defQ}) respectively, and the values
$Q_{FG}(k,\alph,\bet)=Q(0,0,-\infty)$ can be obtained recursively
from (\ref{equ:rec3}) with the distributions of $Y_i$ obtained from
the initial distribution $G$ for $i\in (\alph,\bet]$ and from the
initial distribution $F$ for other values of $i\in\{1,\ldots,k\}$.
Note that the sum of probabilities for all available values $l$ and
$r$ can be larger than $1$.

Let $Y_1,Y_2,\ldots$ be a sequence from the distribution $H$. We use the notations of Section \ref{sec:MLEprop} and note
$$
Q_{H,k}(x,y,s)=\Pr_H(V_k>x,W_k^-\leq y,S_k>s),
$$
$$
Q'_{H,k}(x,y,ds)=\Pr_H(V_k>x,W_k^-\leq y,S_k\in [s,s+ds)),
$$
$$
Q_{H,k}(x,y)=\Pr_H(V_k>x,W_k^-\leq y)=Q_{H,k}(x,y,-\infty),
$$
$$
A_{H,k}(x,dy)=A_{H,k}^+(x,dy)=\Pr_H(R_k^+\leq x,-S_k^+\in [y,y+dy))),
$$
$$
A_{H,k}^-(x,dy)=\Pr_H(R_k^-\leq x,-S_k^-\in [y,y+dy))).
$$
The functions $Q_{H,k}$ and $Q'_{H,k}$ can be obtained recursively by Lemma A.\ref{lem:markov3} under $F_i\equiv H$, $i=1,\ldots,k$.

The probabilities $p_{lr}$ can be obtained from distributions of i.i.d. random variables.

{\it Case 0:} Under $l=0$ and $r=0$,
$$
p_{00}=R_{F,a}(0)Q_{G,b-a}(0,0)R_{F,n-b}(0).
$$

{\it Case 1:} Under $l<b-a+r$ and $r<a-b$,
\arraycolsep=1pt
\begin{eqnarray*}
p_{lr}&=&
\int_0^{\infty}\int_{-x}^{\infty}
\Pr_F(R_{a+l}^+=0)\Pr(V_{d}>0,W^{-}_{d} =0)\\
&& \times\Pr_F(R_{a-b-r}^+=0, -S_{a-b-r}^+\in [x,x+dx)) \\
&& \times\Pr_G(R_{b-a}^+\leq x,-S_{b-a}^+\in [y,y+dy)) \Pr_F(R_{n-b}^+\leq x+y)
\\
&
=&R_{F,a+l}(0)Q_{G,d}(0,0)\\
&\times& \int_0^{\infty}\int_{-x}^{\infty} R_{F,n-b}(x+y)A_{F,a-b-r}(0,dx) A_{G,b-a}(x,dy)
\end{eqnarray*}
where $d=b+r-a-l$.
\smallskip

{\it Case 2:} Under $l<b-a+r$ and $r=a-b$,
\arraycolsep=1pt
\begin{eqnarray*}
p_{lr}&=&
\int_0^{\infty}
\Pr_F(R_{a+l}^+=0)\Pr_F(V_{d}>0,W^{-}_{d} =0)\\
&& \times\Pr_G(R_{b-a}^+=0,-S_{b-a}^+\in [x,x+dx)) \Pr_F(R_{n-b}^+\leq x)
\\
&=&R_{F,a+l}(0)Q_{G,d}(0,0)\int_0^{\infty} R_{F,n-b}(0)(x)A_{G,b-a}(0,dx).
\end{eqnarray*}
\smallskip

{\it Case 3:} Under $l<0$ and $a-b<r<0$,
\arraycolsep=1pt
\begin{eqnarray*}
p_{lr}&=&
\int_0^{\infty}\int_{0}^{\infty}\int_{0}^{\infty}
\Pr_F(R_{a+l}^+=0)\Pr_F(V_{-l}>0,W^{-}_{-l} =0,S_{-l}^+\in [x,x+dx))
\\
&& \times\Pr_G(V_{b+r-a}>-x,W^{-}_{b+r-a} =0,S_{b+r-a}^+\in [y,y+dy))\\
&& \times\Pr_G(R_{-r}^+=0,-S_{-r}^+\in [z,z+dz)) \Pr_F(R_{n-b}^+\leq z)
\\
&
=&R_{F,a+l}(0)\\
&\times& \int_0^{\infty}\int_{0}^{\infty}\int_{0}^{\infty} R_{F,n-b}(z) Q'_{F,-l} (0,y,dx)Q'_{G,b+r-a}(-x,0,dy) A_{G,-r}(0,dz).
\end{eqnarray*}
\smallskip

{\it Case 4:} Under $l<0$ and $r=0$,
\arraycolsep=1pt
\begin{eqnarray*}
p_{lr}&=&
\int_0^{\infty}\int_{0}^{\infty}
\Pr_F(R_{a+l}^+=0)\Pr_F(V_{-l}>0,W^{-}_{-l}\leq y,S_{l}\in [x,x+dx))
\\
&& \times\Pr_G(V_{b-a}>-x,W^-_{b-a}=0,S_{b-a}\in [y,y+dy))\Pr_F(R_{n-b}^+=0)
\\
&
=&R_{F,a+l}(0) R_{F,n-b}(0) \int_0^{\infty}\int_{0}^{\infty} Q'_{F,-l} (0,y,dx)Q'_{G,b-a}(-x,0,dy).
\end{eqnarray*}
and
\begin{eqnarray*}
\Pr_G(V_{b-a}&>&-x,W^-_{b-a}=0,S_{b-a}\in [y,y+dy)) \\
&=&
\Pr_G(V_{b-a}^->0,W_{b-a}\leq x,S_{b-a}^{-}\in [y,y+dy)).
\end{eqnarray*}
\smallskip

{\it Case 5:} Under $l<0$ and $r>0$,
\arraycolsep=1pt
\begin{eqnarray*}
p_{lr}&=&
\int_0^{\infty}\int_{0}^{\infty}\int_{-(x\wedge y)}^{\infty}
\Pr_F(R_{a+l}^+=0)\Pr_F(V_{-l}>0,W^{-}_{-l}\leq y+z,S_{-l}\in [x,x+dx))
\\
&& \times\Pr_G(V_{b-a}>-x,W^{-}_{b-a}\leq y,S_{b-a}\in [z,z+dz))
\\
&& \times\Pr_G(V_{r}>-x-z,W^{-}_{r}=0,S_{r}\in [y,y+dy))\Pr_F(R_{n-b-r}^+=0)
\\
&
=&R_{F,a+l}(0) R_{F,n-b-r}(0)
\\
&& \times\int_0^{\infty}\int_{0}^{\infty}\int_{-(x\wedge y)}^{\infty} \!\!\! Q'_{F,-l} (0,y+z,dx)Q'_{G,r}(-x-z,0,dy)Q'_{G,b-a}(-x,y,dz).
\end{eqnarray*}
\smallskip

{\it Case 6:} Under $l=0$ and $r<0$,
\arraycolsep=1pt
\begin{eqnarray*}
p_{lr}&=&
\int_0^{\infty}
\Pr_F(R_{a}^+=0)\Pr_G(V_{b-a+r}>0,W^-_{b-a+r}=0)
\\
&& \times\Pr_G(R^+_{-r}=0,-S_{-r}\in [y,y+dy))\Pr_F(R_{n-b}^+\leq y)
\\
&=&R_{F,a}(0) Q_{G,d}(0,0)\int_0^{\infty} R_{F,n-b}(x)A_{G,-r} (0,dx).
\end{eqnarray*}
\smallskip

{\it Case 7:} Under $l=0$ and $r>0$,
\arraycolsep=1pt
\begin{eqnarray*}
p_{lr}&=&
\int_0^{\infty}\int_{0}^{\infty}
\Pr_F(R_{a}^+=0)\Pr_G(V_{b-a}>0,W^-_{b-a}\leq y,S_{b-a}\in [x,x+dx))
\\
&& \times\Pr_G(V_{r}>-x,W^-_{r}=0,S_{r}\in [y,y+dy))\Pr_F(R_{n-b}^+=0)
\\
&
=&R_{F,a}(0) R_{F,n-b}(0) \int_0^{\infty}\int_{0}^{\infty} Q'_{G,b-a} (0,y,dx)Q'_{F,r}(-x,0,dy).
\end{eqnarray*}
\smallskip

{\it Case 8:} Under $l>0$ and $r<0$,
\arraycolsep=1pt
\begin{eqnarray*}
p_{lr}&=&
\int_0^{\infty}\int_{0}^{\infty}
\Pr_F(R_{a}^+\leq x)\Pr_G(R_l^+=0,-S_l\in [x,x+dx))
\\
&& \times\Pr_G(V_{b+r-a-l}>0,W^-_{b+r-a-l}=0)\Pr_G(R_r=0,-S_{r}\in [y,y+dy)) \Pr_F(R_{n-b-r}^+\leq y)
\\
&
=&Q_{G,d}(0,0)  \int_0^{\infty} R_{F,a}(x) A_{G,l}(0,dx) \int_{0}^{\infty} R_{F,n-b-r}(y)A_{G,r}(0,dy).
\end{eqnarray*}
\smallskip

{\it Case 9:} Under $l>0$ and $r=0$,
\arraycolsep=1pt
\begin{eqnarray*}
p_{lr}&=&
\int_0^{\infty} \Pr_F(R_{a}^+\leq x)\Pr_G(R_l^+=0,-S_l\in [x,x+dx))
\\
&& \times\Pr_G(V_{b-a-l}>0,W^-_{b-a-l}=0) \Pr_F(R_{F,n-b}^+=0)
\\
&
=&Q_{G,d}(0,0)R_{F,n-b}(0)  \int_0^{\infty} R_{F,a}(x) A_{G,l}(0,dx).
\end{eqnarray*}
\smallskip

{\it Case 10:} Under $0<l<b-a$ and $r>0$,
\arraycolsep=1pt
\begin{eqnarray*}
p_{lr}&=&
\int_0^{\infty}\int_{0}^{\infty}\int_{0}^{\infty}
\Pr_F(R_{a}^+\leq x)\Pr_G(R_l^+=0,-S_l\in [x,x+dx))
\\
&& \times\Pr_G(V_{b-a-l}>0,W^-_{b-a-l}\leq z,S_{b-a-l}\in [y,y+dy))\\
&& \times\Pr_G(V_{r}>-y,W^-_{r}=0,S_{r}\in [z,z+dz)) \Pr_F(R_{n-b-r}^+\leq 0)
\\
&
=& R_{F,n-b-r}(0) \int_0^{\infty} R_{F,a}(x) A_{G,l}(0,dx) \int_0^{\infty}\int_0^{\infty} Q'_{G,b-a-l}(0,z,dy) Q'_{G,r}(-y,0,dz) .
\end{eqnarray*}
\smallskip

{\it Case 11:} Under $l=b-a$ and $r>0$,
\arraycolsep=1pt
\begin{eqnarray*}
p_{lr}&=&
\int_0^{\infty}
\Pr_F(R_{a}^+\leq x)\Pr_F(R_{b-a}^+=0,-S_{b-a}\in [x,x+dx))
\\
&& \times\Pr_G(V_{r}>0,W^-_{r}\leq 0)\Pr_F(R_{n-b-r}^+\leq 0)
\\
&
=& R_{F,n-b-r}(0)Q_{G,r}(0,0) \int_0^{\infty} R_{F,a}(x) A_{G,l}(0,dx).
\end{eqnarray*}
\smallskip

{\it Case 12:} Under $l>b-a$ and $r>0$,
\arraycolsep=1pt
\begin{eqnarray*}
p_{lr}&=&
\int_0^{\infty}\int_{-x}^{\infty}
\Pr_F(R_{a}^+\leq x+y)\Pr_G(R_{b-a}^+\leq x,-S_{b-a}\in [y,y+dy))
\\
&& \times\Pr_G(R_{a-b+l}^+=0,-S_{a-b+l}\in [x,x+dx))\\
&& \times\Pr_F(V_{d}>0,W^-_{d}=0) \Pr_F(R_{n-b-r}^+\leq 0)
\\
&
=& R_{F,n-b-r}(0)Q_{F,d}(0,0) \int_0^{\infty}\int_{-x}^{\infty} R^+_{a}(x+y)A_{G,a-b+l}(0,dx)A_{F,b-a}(x,dy).
\end{eqnarray*}



\begin{thebibliography}{10}

\bibitem{Baron00a}
M.~Baron.
\newblock Nonparametric adaptive change-point estimation and on-line detection.
\newblock {\em Sequential Analysis}, 19(12):1--23, 2000.

\bibitem{BaronRosenbergSidorenko01}
M.~Baron, M.~Rosenberg, and N.~Sidorenko.
\newblock Electricity pricing: modeling and prediction with automatic spike
  detection.
\newblock {\em Energy, Power, and Risk Management}, October 2001:36--39, 2001.

\bibitem{BaronRosenbergSidorenko02}
M.~Baron, M.~Rosenberg, and N.~Sidorenko.
\newblock Divide and conquer: forecasting power via automatic price regime
  separation.
\newblock {\em Energy, Power, and Risk Management}, Match 2002:70--73, 2002.

\bibitem{BaronRukhin97}
M.~Baron and A.~L. Rukhin.
\newblock Asymptotic behavior of confidence regions in the change-point
  problem.
\newblock {\em J. of Stat. Planning and Inference}, 58:263--282, 1997.

\bibitem{bianchi1993time}
A.~M. Bianchi, L.~Mainardi, E.~Petrucci, M.~G. Signorini,
M.~Mainardi, and
  S.~Cerutti.
\newblock Time-variant power spectrum analysis for the detection of transient
  episodes in hrv signal.
\newblock {\em IEEE Transactions on Biomedical Engineering}, 40(2):136--144,
  1993.

\bibitem{crosier1986new}
Ronald~B Crosier.
\newblock A new two-sided cumulative sum quality control scheme.
\newblock {\em Technometrics}, 28(3):187--194, 1986.

\bibitem{egea2018performance}
D.~Egea-Roca, J.~A. L{\'o}pez-Salcedo, G.~Seco-Granados, and H.~V.
Poor.
\newblock Performance bounds for finite moving average tests in transient
  change detection.
\newblock {\em IEEE Transactions on Signal Processing}, 66(6):1594--1606, 2018.

\bibitem{guepie2017detecting}
B.~K. Gu{\'e}pi{\'e}, L.~Fillatre, and I.~Nikiforov.
\newblock Detecting a suddenly arriving dynamic profile of finite duration.
\newblock {\em IEEE Transactions on Information Theory}, 63(5):3039--3052,
  2017.

\bibitem{guepie2012sequential}
B.~K. Gu{\'e}pi{\'e}, L.~Fillatre, and I.~V. Nikiforov.
\newblock Sequential detection of transient changes.
\newblock {\em Sequential Analysis}, 31(4):528--547, 2012.

\bibitem{Hinkley70}
D.~V. Hinkley.
\newblock Inference about the change-point in a sequence of random variables.
\newblock {\em Biometrika}, 57:1--17, 1970.

\bibitem{HochbergTamhane87}
Y.~Hochberg and A.~C. Tamhane.
\newblock {\em Multiple comparison procedures}.
\newblock Wiley, New York, 1987.

\bibitem{HuRukhin95}
I.~Hu and A.~L. Rukhin.
\newblock A lower bound for error probability in change-point estimation.
\newblock {\em Statistica Sinica}, 5:319--331, 1995.

\bibitem{JJSiegmund92}
B.~James, K.~L. James, and D.~Siegmund.
\newblock Asymptotic approximations for likelihood ratio tests and confidence
  regions for a change-point in the mean of a multivariate normal distribution.
\newblock {\em Statist. Sinica}, 2:69--90, 1992.

\bibitem{Leadbetter83}
M.~R. Leadbetter, G.~Lindgren, and H.~Rootzen.
\newblock {\em Extremes and Related Properties of Random Sequences and
  Processes}.
\newblock Springer-Verlag, New York, 1983.

\bibitem{noonan2020power}
J.~Noonan and A.~Zhigljavsky.
\newblock Power of the mosum test for online detection of a transient change in
  mean.
\newblock {\em Sequential Analysis}, 39(2):269--293, 2020.

\bibitem{Page54}
E.~S. Page.
\newblock Continuous inspection schemes.
\newblock {\em Biomterika}, 41:100--115, 1954.

\bibitem{repin1991detection}
V.~G. Repin.
\newblock Detection of a signal with unknown moments of appearance and
  disappearance.
\newblock {\em Problemy Peredachi Informatsii}, 27(1):61--72, 1991.

\bibitem{reynolds2010evaluation}
J.~Reynolds, R.~Marion, and J.~Lou.
\newblock An evaluation of a glr control chart for monitoring the process mean.
\newblock {\em Journal of quality technology}, 42(3):287--310, 2010.

\bibitem{RosenbergBryngelsonBaronPapalexopoulos10}
M.~Rosenberg, J.~D. Bryngelson, M.~Baron, and A.~D. Papalexopoulos.
\newblock Transmission valuation analysis based on real options with price
  spikes.
\newblock In {\em S. Rebennack, P.M. Pardalos, M.V.F. Pereira and N. Iliadis,
  eds. Handbook of Power Systems II; Energy Systems Part I}, pages 101--125,
  Springer, Berlin-Heiderberg, 2010.

\bibitem{RosenbergBryngelsonSidorenkoBaron02}
M.~Rosenberg, J.~D. Bryngelson, N.~Sidorenko, and M.~Baron.
\newblock Price spikes and real options: transmission valuation.
\newblock In {\em E. I. Ronn, ed., Real Options and Energy Management}, pages
  323--370, Risk Books, London, 2002.

\bibitem{Shiryaev95}
A.~N. Shiryaev.
\newblock {\em Probability, 2nd edition}.
\newblock Springer-Verlag, New York, 1995.

\bibitem{stroock2013mathematics}
D.~W. Stroock.
\newblock {\em Mathematics of probability}, volume 149.
\newblock American Mathematical Soc., 2013.

\bibitem{tafakori2018forecasting}
L.~Tafakori, A.~Pourkhanali, and F.~A. Fard.
\newblock Forecasting spikes in electricity return innovations.
\newblock {\em Energy}, 150:508--526, 2018.

\bibitem{tartakovskii1988detection}
A.~G. Tartakovskii.
\newblock Detection of signals with random moments of appearance and
  disappearance.
\newblock {\em Problemy Peredachi Informatsii}, 24(2):39--50, 1988.

\bibitem{tartakovsky2021optimal}
A.~G. Tartakovsky, N.~R. Berenkov, A.~E. Kolessa, and I.~V.
Nikiforov.
\newblock Optimal sequential detection of signals with unknown appearance and
  disappearance points in time.
\newblock {\em IEEE Transactions on Signal Processing}, 69:2653--2662, 2021.

\bibitem{Woodroofe82}
M.~Woodroofe.
\newblock {\em Nonlinear renewal theory in sequential analysis}.
\newblock SIAM, 1982.

\bibitem{zhang2017regime}
L.~Zhang and Y.~Li.
\newblock Regime-switching based vehicle-to-building operation against
  electricity price spikes.
\newblock {\em Energy Economics}, 66:1--8, 2017.

\bibitem{zhou2019improving}
B.~Zhou, M.~Chioua, M.~Bauer, J.~C. Schlake, and N.~F. Thornhill.
\newblock Improving root cause analysis by detecting and removing transient
  changes in oscillatory time series with application to a 1, 3-butadiene
  process.
\newblock {\em Industrial \& Engineering Chemistry Research},
  58(26):11234--11250, 2019.

\end{thebibliography}


\end{document}